\documentclass{amsart}

\newtheorem{definition}{Definition}[section]
\newtheorem{proposition}[definition]{Proposition}

\title{Remarks on the Definition of a Courant Algebroid}

\begin{document}
\maketitle

\begin{center}
{\large 
Kyousuke Uchino \\
Department of Mathematics\\
Tokyo University of Science (Science University of Tokyo)\\
Kagurazaka 1-3, Shinjuku-ku, Tokyo, 162-8601, Japan}
\end{center}

\begin{abstract}
The notion of Courant algebroid was introduced by Liu, Weinstein and Xu
in 1997. Its definition consists of five axioms and an defining relation for a 
derivation. It is shown that two of the axioms and the relation (assuming only the Leibniz rule) follow from the rest of the axioms.
\end{abstract}

\section{Introduction}

The Courant bracket on $TM \oplus T^{*}M$, where $M$ is a smooth manifold,
was defined in 1990 by T.Courant \cite{Cou}. A similar bracket was defined on the double of any Lie bialgebroid by Liu, Weinstein, and Xu \cite{Liu}. They show that the double is not a Lie algebroid, but a more complicated object which they call a Courant algebroid.\\
\\
{\large {\bf Definition of a Courant algebroid}} \cite{Liu} \\
A Courant algebroid is a vector bundle $E \rightarrow M$ equipped with a 
nondegenerate symmetric bilinear form $(\hspace{0.5em},\hspace{0.5em})$, 
a skew-symmetric bracket $[\hspace{0.5em},\hspace{0.5em}]$ on $\Gamma E$ and 
a bundle map $\rho :E \rightarrow TM$ satisfying the following relations:\\
\begin{enumerate}
\item[[C1]]\ $[[x , y] , z]+[[z , x], y]+[[y, z] , x]=DT(x, \ y, \ z) \quad
\forall x,y,z \in \Gamma E,$
\\
\item[[C2]]\ $\rho [x, y]=[\rho(x), \rho(y)] \quad \forall x,y \in \Gamma E,$
\\
\item[[C3]]\ $[x , fy]=f[x , y]+(\rho(x)f)y-(x, \ y)Df \quad 
\forall x,y \in \Gamma E, \forall f \in C^{\infty}(M),$
\\
\item[[C4]]\ $\rho \circ D=0,$ i.e., $(Df,\ Dg)=0 \quad
\forall f,g \in C^{\infty}(M),$
\\
\item[[C5]]\ $\rho(x)(y,\ z)=([x, y]+D(x,\ y),\ z)+(y,\ [x, z]+D(x,\ z))\ 
\forall x,y,z \in \Gamma E,$
\end{enumerate}
\quad \\
where $T(x,\ y,\ z)=\frac{1}{3}(([x , y], z)+([y , z], x)+([z, x], y))$,
and $D:C^{\infty}(M) \rightarrow \Gamma E$ is the map given by the composition of the dual map $\rho^{*}$ of $\rho$ and the usual differential $d$, 
$D=\frac{1}{2}\beta^{-1}\rho^{*}d$, $\beta$ being the isomorphism between 
$E$ and $E^{*}$ given by the bilinear form. Here $\Gamma E$ denotes the 
space of all smooth sections of $E$.\\

\noindent\textit{Remark 1}. $D=\frac{1}{2}\beta^{-1}\rho^{*}d$ is equivalent to 
$(Df,\ x)=\frac{1}{2}\rho(x)f \ \forall x \in \Gamma E, 
\ \forall f \in C^{\infty}(M)$.\\

In this paper we show that [C3] and [C4] and the property of $D$, i.e.,
$D=\frac{1}{2}\beta^{-1}\rho^{*}d$, follow from [C2], [C5] and the Leibniz 
rule for $D$, (L) below.

\section{Main Results}
We assume [C2] and [C5] in the definition above, and only that $D$ 
is a map from $C^{\infty}(M)$ to $\Gamma E$ satisfying [C5] and the 
Leibniz rule:
\begin{equation}
D(fg)=fD(g)+gD(f),  \quad {\forall}f,g \in C^{\infty}(M). 
\tag{L}
\end{equation}
We then have
\begin{proposition}
$$(i) \quad [x , fy]=f[x , y]+\rho(x)(f)y-(x, \ y)Df,\
\forall x,y \in \Gamma E, \forall f \in C^{\infty}(M),$$ 
$$(ii) \quad \rho \circ D=0, \forall x,y \in \Gamma E.$$ 
\begin{proof}
First we show {\em (i)}. By [C5] we have
\begin{eqnarray*}
\rho(x)(fy, \ z)&=&([x , fy]+D(x \ fy),\ z)+(fy, \ [x , z]+D(x \ z)).
\end{eqnarray*}
On the other hand, the Leibniz rule gives 
\begin{eqnarray*}
\rho(x)(fy, \ z)&=&f \rho(x)(y \ z)+(y \ z)\rho(x)(f),\\
   &=&(f[x , y]+fD(x \ y), \ z)+(fy, \ [x , z]+D(x \ z))+(y \ z)\rho(x)(f).
\end{eqnarray*}
Hence we have 
\begin{eqnarray*}
(f[x , y]+fD(x, \ y), \ z)+(y, \ z)\rho(x)(f)&=&([x , fy]+D(x, \ fy), \ z).
\end{eqnarray*}
Since (\hspace{0.5em},\hspace{0.5em}) is nondegenerate we obtain 
\begin{eqnarray*}
f[x , y]+fD(x, \ y)+\rho(x)(f)y&=&[x , fy]+D(x, \ fy).
\end{eqnarray*}
By the Leibniz rule we have $D(x, \ fy)=fD(x, \ y)+(x, \ y)Df$
and then get {\em (i)}, $[x , fy]=f[x , y]+\rho(x)(f)y-(x, \ y)Df.$\\
\\
Next we show {\em (ii)}. By [C2] we have 
$$\rho [x , fy]\ =\ [\rho(x) , f\rho(y)]\ =\ 
f[\rho(x) , \rho(y)]+\rho(x)(f)\rho(y).$$
Then applying $\rho$ to the both sides of {\em (i)} gives the desired result.
\end{proof}
\end{proposition}
Next we show a map $D$ satisfying [C5] and the Leibniz rule (L) can be 
written $D=\frac{1}{2}\beta^{-1}\rho^{*}d$. 
By Remark 1, it suffices to show that
\begin{proposition}
\[
(Df, \ y)=\frac{1}{2}\rho(y)(f), \quad \forall y \in \Gamma E, 
\ \forall f \in C^{\infty}(M).
\]
\end{proposition}
\begin{proof}
If $y=0$ the identity is trivial. We assume $y \neq 0$.\\

$\forall x,y \in \Gamma E \ \forall f \in C^{\infty}(M)$,  [C5] gives
\begin{eqnarray*}
\rho(fx)(y, \ y)&=&([fx , y]+D(fx, \ y), \ y)+([fx , y]+D(fx, \ y), \ y).
\end{eqnarray*}
We apply Proposition 1 {\em (i)} and the Leibniz rule for $D$ to the line
above and get
\begin{eqnarray*}
\rho(fx)(y, \ y)&=&(f[x , y]-\rho(y)(f)x+(x,\ y)Df+fD(x,\ y)+(x,\ y)Df,\ y)\\ 
         &+&(f[x , y]-\rho(y)(f)x+(x, \ y)Df+fD(x, \ y)+(x, \ y)Df, \ y),\\
         &=&f([x , y]+D(x, \ y), \ y)+f([x , y]+D(x, \ y), \ y)\\
         &+&(-\rho(y)(f)x+2(x, \ y)Df, \ y )+(-\rho(y)(f)x)+2(x, \ y)Df, \ y).
\end{eqnarray*}
Again by [C5]
\begin{eqnarray*}
\rho(fx)(y, \ y)&=&f\rho(x)(y, \ y)+2(-\rho(y)(f)x+2(x, \ y)Df, \ y ).
\end{eqnarray*}
Hence we have
$0=(-\rho(y)(f)x+2(x, \ y)Df, \ y )\ = \ -\rho(y)(f)(x, \ y)+2(x,\ y)(Df,\ y)$.
Since (\hspace{0.5em},\hspace{0.5em}) is nondegenerate we obtain 
$$ 0=-\rho(y)(f)y+2(Df, \ y)y,\ =\ (2(Df, \ y)-\rho(y)(f))y,$$
which gives $(Df, \ y)=\frac{1}{2}\rho(y)(f)$.
\end{proof}
\noindent \textit{Remark 2.1. \
Conversely, the Leibniz rule (L) can be derived from Proposition 1(i)}.
Namely, if a map $D:C^{\infty}(M) \rightarrow \Gamma E$ satisfies the 
identity in Proposition 1(i), then $D$ fulfills the Leibniz rule (L).\\

In fact, by Proposition 1 {\em (i)} we have
\begin{eqnarray*}
(x,\ y)Df=f[x , y]+(\rho(x)f)y-[x , fy].
\end{eqnarray*}
Since right-hand side is in $\Gamma E$ and $(x,\ y)$  in 
$C^{\infty}(M)$, we obtain $Df$ $\in \Gamma E$. Applying 
Proposition 1{\em (i)} to $[x , (fg)y]$ and $[x , f(gy)]$, 
we obtain $D(fg)=fD(g)+gD(f)$.\\
\\
\noindent \textit{Remark 2.2.}
For a Lie algebroid, Grabowski and Marmo showed that a homomorphism property of $\rho$ with respect to [\, ,\, ], i.e., [C2], follows from the rest of the axioms of a Lie algebroid(\cite{Gra}). It is thus natural to ask if, in a Courant algebroid, one can derive [C2] from the rest of the axioms. At present, this is still uncertain. We can however show the following.\\
\quad If we assume only [C5] and (L), then
\begin{equation}
\rho(x)\rho(y)(f)=\rho(x \circ y)(f)+\rho(y)\rho(x)(f)+2([x , Df]-D(x, \ Df), \ y),
\tag{A}
\end{equation}
where $x \circ y:=[x , y]+D(x, \ y)$. Indeed by [C5] and Proposition 2, we have
\begin{eqnarray*}
\rho(x)\rho(y)(f)&=&2(x \circ y,\ Df)+2(y,\ x \circ Df)+\rho(y)\rho(x)(f)-\rho(y)\rho(x)(f).
\end{eqnarray*}
Again by Proposition 2, we obtain $2(x \circ y,\ Df)=\rho(x \circ y)(f)$ and $\rho(y)\rho(x)(f)=4(y,\ D(x,\ Df))$. Then we have (A).\\
\quad {\it By (A) we see [C2] is equivalent to [x , Df]=D(x,\ Df) as follows.}
Since $[x , y]=\frac{1}{2}(x \circ y -y \circ x)$, $\rho [x , y]=\frac{1}{2}(\rho(x \circ y) -\rho(y \circ x))$. If $[x , Df]=D(x,\ Df)$ then (A) gives $\rho(x \circ y)=[\rho(x) , \rho(y)]$ and then we have [C2]. Conversely, if [C2] is assumed then $\rho \circ D=0$ by Proposition 1 $(ii)$, which gives $\rho(x \circ y)=\rho([x \ y]+D(x,\ y))=\rho [x , y]=[\rho(x) , \rho(x)]$. Hence $[x , Df]=D(x, Df).$\\
\quad It is to be noticed that a Courant algebroid has the property $[x , Df]=D(x, Df)$ (\cite{Roy}, \cite{Wei}).\\
\\
\noindent \textit{Remark 2.3.}
It is  known that a Courant algebroid is given equivalently by the following axioms.\\
\\
{\large {\bf Alternative definition of a Courant algebroid}} \cite{Roy} \\
A Courant algebroid is a vector bundle $E \rightarrow M$ equipped with a 
nondegenerate symmetric bilinear form $(\hspace{0.5em},\hspace{0.5em})$, 
a non-skew-symmetric bilinear operation $\circ$ on $\Gamma E$ and 
a bundle map $\rho :E \rightarrow TM$ satisfying the following relations:\\
\begin{enumerate}
\item[[C'1]]\ $x \circ (y \circ z)=(x \circ y) \circ z + y \circ (x \circ z)
\quad
\forall x,y,z \in \Gamma E,$
\\
\item[[C'2]]\ $\rho (x \circ y)=[\rho(x), \rho(y)] \quad \forall x,y \in \Gamma E,$
\\
\item[[C'3]]\ $x \circ fy=f(x \circ y)+(\rho(x)f)y \quad 
\forall x,y \in \Gamma E, \forall f \in C^{\infty}(M),$
\\
\item[[C'4]]\ $x \circ x=D(x,\ x) \quad \forall x \in \Gamma E,$
\\
\item[[C'5]]\ $\rho(x)(y,\ z)=(x \circ y,\ z)+(y,\ x \circ z)\ 
\forall x,y,z \in \Gamma E,$
\end{enumerate}
\quad \\
where $D:C^{\infty}(M) \rightarrow \Gamma E$ is the map defined by 
$D=\frac{1}{2}\beta^{-1}\rho^{*}d$, $\beta$ being the isomorphism between 
$E$ and $E^{*}$ given by the bilinear form. Here $\Gamma E$ denotes the 
space of all smooth sections of $E$.\\
\\
{\it The above axioms [C'2], [C'3] and assumption for the derivation $D=\frac{1}{2}\beta^{-1}\rho^{*}d$ are also derived from the rest of the axioms and condition (L) in a similar way.}\\
\begin{center}
Acknowledgements\\
\end{center}
\quad I am grateful to the referee for a helpful comment and I would also like to thank very much Professor Daniel Sternheimer and Professor Akira Yoshioka for encouragement.


\end{document}